\documentclass{amsart}

\usepackage{amssymb}
\usepackage[xdvi]{graphics}

\theoremstyle{plain}
\newtheorem{thm}{Theorem}
\newtheorem{prop}[thm]{Proposition}
\newtheorem{lemma}[thm]{Lemma}
\newtheorem{cor}[thm]{Corollary}
\newtheorem*{claim}{Claim}

\theoremstyle{definition}
\newtheorem{defn}[thm]{Definition}

\theoremstyle{remark}
\newtheorem{rmk}[thm]{Remark}

\newcommand{\eps}{\varepsilon}

\newcommand{\C}{\mathbb{C}}
\newcommand{\R}{\mathbb{R}}

\newcommand{\N}{\mathbb{N}}
\newcommand{\Z}{\mathbb{Z}}
\newcommand{\T}{\mathcal{T}}
\newcommand{\A}{\mathcal{A}}
\newcommand{\B}{\mathcal{B}}
\newcommand{\K}{\mathfrak{K}}
\newcommand{\SR}{\Sigma_R}
\renewcommand{\Re}{\operatorname{Re}}
\newcommand{\act}{\diamond}
\newcommand{\lact}{\cdot}
\newcommand{\altbullet}{-}
\renewcommand{\u}{\mathcal{U}}

\DeclareMathOperator{\Op}{Op}
\DeclareMathOperator{\Fin}{Fin}

\begin{document}

\title[Property A, Nuclearity, and Uniform Embeddings]{Property A, partial translation structures \\ and Uniform Embeddings in Groups}

\author{J. Brodzki}
\author{G.A. Niblo}
\author{N.J. Wright}

\address{School of Mathematics, University of Southampton, Southampton SO17 1BJ}


\begin{abstract}
We define the  concept of a  partial   translation structure $\T$ on a metric space  $X$ and we 
 show that there is a natural $C^*$-algebra $C^*(\T)$ associated with it  which is a subalgebra of the uniform Roe algebra $C^*_u(X)$. We introduce a coarse invariant of the metric which provides an obstruction to embedding  the space in a group.  When the space is sufficiently group-like, as determined by our invariant, properties of the Roe algebra can be deduced from those of $C^*(\T)$. We also give a proof
of the fact  that the uniform Roe algebra of a metric space is a coarse invariant up to 
Morita equivalence. 
\end{abstract}

\maketitle

Many interesting geometric properties of spaces and groups are captured by the structure of $C^*$-algebras associated with those objects. For example, a discrete group $G$ is amenable if and only if 
the full $C^*$-algebra $C^*(G)$ is nuclear  \cite{Lance}.

In a similar vein, for a discrete group $G$, Yu's property A is equivalent  both to  the nuclearity of the 
uniform Roe algebra $C^*_u(G)$ and to the exactness of the reduced $C^*$-algebra 
$C^*_r(G)$. This follows from the results of Anantharaman-Delaroche and Renault \cite{AD}, Higson and Roe \cite{HR}, Guentner and  Kaminker \cite{GK}, and Ozawa \cite{Oz}. 
While property A and the uniform Roe algebra can be defined for arbitrary metric spaces, we cannot 
generalise these results without a good analogue of the reduced $C^*$-algebra
of a group. 

In this paper we introduce a $C^*$-algebra to fulfill this role. To do so we carry out the following programme. First we define the notion of 
a partial translation structure (Definition \ref{pt}) on a uniformly discrete metric space, which captures geometrically 
the interplay between the left and the right action of a group on itself. In broad terms, this can be 
described as follows. 
In  Euclidean space translations are distinguished from other isometries of the space by the fact that they move each point by the same distance. Let us assume that a group $G$ is equipped with a left invariant metric $d$, which means that for any elements $g, s,t$ of $G$,
$d(gs,gt) = d(s,t)$. In other words, multiplication on the left acts by isometries on the metric space 
$(G, d)$. On the other hand, the right multiplication by a fixed element $g$ of $G$ moves each element $r$ of $G$ by the same distance: 
$
d(r,rg) = d(e,g)
$, where $e$ is the identity of $G$. 
Thus we say that the right multiplication acts by translations on $(G,d)$ even though these translations are not isometries in general.  However, there is an interesting connection between the translations arising from the right 
multiplication and the isometries given by the left multiplication. This interaction is encoded 
in our definition of partial translation structure. It will follow directly from this definition that any group
admits a canonical partial translation structure given by the left and right multiplication. This partial translation structure
is in some sense the best possible, as discussed  in Section \ref{Section3}. To measure how group-like is  a given metric space $X$,  we  introduce a coarse  invariant, a real valued function 
$\kappa_X(R)$, associated to  the space $X$.  When 
$\kappa_X(R)=1$ for all $R>0$ we say that the partial translation structure is free. 

Next we construct 
a $C^*$-algebra, the partial translation algebra $C^*(\T)$, associated to any partial translation structure
$\T$. This is a subalgebra of the uniform Roe algebra $C^*_u(X)$. 
We show that when  the space $X$ is a group and $\T$ is the 
canonical partial translation structure then  $C^*(\T)$ is isomorphic to the reduced $C^*$-algebra of 
the group. We study the analytic properties of the partial translation algebras and we prove 
that when the space admits a free, globally controlled partial translation structure $\T$ then the following are equivalent (Theorem \ref{main}): 
\begin{enumerate}
\item $X$ has property $A$; 
\item $C^*(\T)$ is exact; 
\item $C^*_u(X)$ is exact; 
\item $C^*_u(X)$ is nuclear. 
\end{enumerate}
This is related to a theorem of Skandalis, Tu and Yu, \cite{STY}, who proved the equivalence of 1 and 4 for bounded geometry metric spaces.

In Theorem \ref{ptforgroups} we prove the following. When the space $X$ admits an injective uniform embedding into a group, one can pull back the 
canonical partial translation structure from the group to the space to obtain a free partial translation structure on $X$ which satisfies the hypotheses of the above Theorem \ref{main}.  We conclude furthermore that 
$\kappa_X(R)=1$ for all $R$ in this case. It follows that if $\kappa_X(R) > 1$ for some $R$, then $X$ does not admit an injective uniform embedding into any countable group. Our invariant therefore provides an
obstruction to the existence of such an embedding (Corollary \ref{cor20}). 

The question of the existence of a uniform embedding in a group can be treated locally.  We introduce the notion of a local embedding in a group and show that a discrete metric space $X$ is embeddable in a
group if and only if it is locally embeddable (Theorem \ref{localue}). The main interest in this notion comes from the existence of a 
space $\Gamma_\u$ with the universal property that every bounded geometry space 
locally uniformly embeds in $\Gamma_\u$. We are therefore able to show that $\Gamma_\u$
embeds uniformly in a countable discrete group $G$ if and only if 
for each  uniformly discrete bounded geometry metric space $X$ there 
is a  countable discrete group $G_X$ such that $X$ uniformly embeds in $G_X$.

We remark that property A of a metric space $X$ is a coarse invariant, but the uniform Roe algebra $C^*_u(X)$ is not. We prove, however, that if $X$ and $Y$ are coarsely equivalent then 
$C^*_u(X)$ and $C^*_u(Y)$ are Morita equivalent (Theorem \ref{moritainvariance}). Thus exactness and nuclearity of the uniform 
Roe algebra are preserved by coarse equivalence.

\section{Characterisations of Property A}

In this section we will review various characterisations of property A for metric spaces. The original definition, due to Yu, is stated as follows. 
\begin{defn}[Yu, \cite{Yu}]
A uniformly discrete metric space $(X,d)$ has \emph{property $A$} if for all $R,\eps>0$ there exists a family of finite non-empty subsets $A_x$ of $X\times \N$, indexed by $x$ in $X$, such that
\begin{itemize}
\renewcommand{\labelitemi}{\altbullet}
\item for all $x,y$ with $d(x,y)<R$ we have $\frac{|A_x\Delta A_{y}|}{|A_x\cap A_{y}|}<\eps$;
 
\item there exists $S$ such that for all $x$ and $(y,n) \in A_x$ we have $d(x,y)\leq S$.
\end{itemize}
\end{defn}

We will begin by giving some conditions which are equivalent to property A. For this the following 
terminology will be useful. 

\begin{defn}
A function from $X$ to a Banach space $x\mapsto \xi_x$ has \emph{$(R,\eps)$-variation} if $d(x,y)\leq R$ implies $\|\xi_x-\xi_{y}\|<\eps$.

A kernel $u\colon X\times X \to \R$ or $\C$ has \emph{$(R,\eps)$-variation} if $d(x,y)\leq R$ implies $|u(x,y)-1|<\eps$.

A kernel $u\colon X\times X \to \R$ or $\C$ has \emph{finite propagation} if there exists $R\geq 0$ such that $u(x,y)=0$ for $d(x,y)>R$. The $\emph{propagation}$ of $u$ is the smallest such $R$.
\end{defn}

 We recall that a uniformly discrete metric space $X$ has \emph{bounded geometry} if for all $R>0$ there exists $N$, such that the cardinality of $B_R(x)$ is at most $N$ for all $x$ in $X$.
Hence, if  $u\colon X\times X \to \C$ is a finite propagation kernel on a bounded geometry discrete metric space then there exists $N$ such that  for each $x$ there are no more than $N$ points $y$ such that
 $u(x,y)\neq 0$. Thus $u$ defines a bounded linear map from $l^2(X)$ to itself, $(u*\xi)(x)=\sum_{y\in X} u(x,y)\xi(y)$. These linear maps are also said to have finite propagation. 

The \emph{uniform Roe algebra}, $C^*_u(X)$, is the $C^*$-algebra completion of the algebra of bounded operators on $l^2(X)$ having finite propagation.

The following result summarises 
various characterisations of Property A.

\begin{thm}
\label{propertyAcharacterisations}
Let $X$ be a uniformly discrete bounded geometry  metric space. The following are equivalent.
\begin{enumerate}
\item $X$ has property A.
\label{propA}

\item For all $R,\eps$ there exists a family of vectors $\xi_x \in l^1(X)$ for $x\in X$, such that: 
\begin{itemize}
\renewcommand{\labelitemi}{\altbullet}
\item $\|\xi_x\|_1=1$; the family $(\xi_x)$ has $(R,\eps)$-variation; 
\item there exists $S$ such that for all $x$, $\xi_x$ is supported in the $S$-ball about $x$.
\end{itemize}
\label{l1}

\item For all $R,\eps$ there exist vectors $\xi_x \in l^2(X)$ for $x\in X$, such that: 
\begin{itemize}
\renewcommand{\labelitemi}{\altbullet}
\item $\|\xi_x\|_2=1$; the family $(\xi_x)$ has $(R,\eps)$-variation; 
\item there exists $S$ such that for all $x$, $\xi_x$ is supported in the $S$-ball about $x$.
\end{itemize}
\label{l2}

\item For all $R,\eps$ there exist vectors $\xi_x \in l^2(X)$ for $x\in X$, such that: 
\begin{itemize}
\renewcommand{\labelitemi}{\altbullet}
\item $\|\xi_x\|_2=1$; the family $(\xi_x)$ has $(R,\eps)$-variation; 
\item for all $\delta>0$ there exists $S$ such that for all $x$, the restriction of $\xi_x$ to the $S$-ball about $x$ has norm at least $1-\delta$.
\end{itemize}
\label{l2delta}

\item There exists $\delta<1$ such that for all $R,\eps$ there exist vectors $\xi_x \in l^2(X)$ for $x\in X$, such that: 
\begin{itemize}
\renewcommand{\labelitemi}{\altbullet}
\item $\|\xi_x\|_2=1$; the family $(\xi_x)$ has $(R,\eps)$-variation; 
\item there exists $S$ such that for all $x$ the restriction of $\xi_x$ to the $S$-ball about $x$ has norm at least $1-\delta$, and the restriction of $\xi_x$ to the set $B_{R+S}(x)\setminus B_S(x)$ has norm at most $\eps$.
\end{itemize}
\label{l2deltaweak}

\item For all $R,\eps$ there exists a Hilbert space $H$ and vectors $\xi_x \in H$ for $x\in X$, such that:
\begin{itemize}
\renewcommand{\labelitemi}{\altbullet}
\item  $\|\xi_x\|_2=1$; the family $(\xi_x)$ has $(R,\eps)$-variation; 
\item there exists $S$ such that $d(x,y)>S$ implies $\langle\xi_x,\xi_{y}\rangle=0$.
\end{itemize}
\label{Hilbert}

\item For all $R,\eps$ there exists a finite propagation positive type kernel $u\colon X\times X \to \R$ such that  $u$ has $(R,\eps)$-variation.
\label{Ozawakernel}

\item For all $R,\eps$ there exists a positive type kernel $u\colon X\times X \to \C$ such that: $u$ has $(R,\eps)$-variation; convolution with $u$ defines a bounded operator in the uniform Roe algebra $C^*_u(X)$.
\label{Roekernel}
\end{enumerate}
\end{thm}

\begin{proof}
The equivalence of \ref{propA} and \ref{l1} is proved in \cite{HR}, Lemma 3.5. The equivalence of \ref{l1} and \ref{l2} is proved in \cite{Tu} Proposition 3.2.

\ref{l2} $\iff$ \ref{l2delta} $\iff$ \ref{l2deltaweak}: Note that \ref{l2} $\implies$ \ref{l2delta} trivially. Given \ref{l2delta}, fix $\delta<1$, and for any $R,\eps$, take $\delta'$ to be the smaller of $\delta,\eps$. For $S$, $\xi_x$ such that the restriction of $\xi_x$ to the $S$-ball has norm at least $1-\delta'\geq 1-\delta$, the restriction to $B_{R+S}(x)\setminus B_S(x)$ has norm at most $\delta'\leq \eps$. Thus \ref{l2delta} $\implies$ \ref{l2deltaweak}.

To show \ref{l2deltaweak} $\implies$ \ref{l2}, fix $\delta$ and given any $R,\eps$ let $\xi_x$ and $S$ be as in \ref{l2deltaweak}. Define $\zeta_x$ to be the restriction of $\xi_x$ to the $R+S$-ball about $x$. We write $\zeta_x-\zeta_y$ as the sum of three parts: $\xi_x$ restricted to $B_{R+S}(x)\setminus B_{R+S}(y)$; $\xi_x-\xi_y$ restricted to $B_{R+S}(x)\cap B_{R+S}(y)$; and $-\xi_y$ restricted to $B_{R+S}(y)\setminus B_{R+S}(x)$. Note that for $d(x,y)<R$ we have $B_{R+S}(x)\setminus B_{R+S}(y)$ contained in $B_{R+S}(x)\setminus B_S(x)$, thus the restriction of $\xi_x$ to this set has norm at most $\eps$. Similarly for the corresponding restriction of $-\xi_y$. The variation condition ensures that any restriction of $\xi_x-\xi_y$ has norm at most $\eps$, hence for $d(x,y)<R$ we have $\|\zeta_x-\zeta_{y}\|\leq 3\eps$. Let $\eta_x=\zeta_x/\|\zeta_x\|$. Now use the estimates
\begin{align*}
\| \eta_x - \eta_{y} \| &\leq \frac 1{\|\zeta_x\|}\|\zeta_x-\zeta_{y}\| + \|\zeta_{y}\|\left|\frac 1{\|\zeta_x\|}-\frac 1{\|\zeta_{y}\|}\right| \\
&= \frac 1{\|\zeta_x\|}\|\zeta_x-\zeta_{y}\| + \frac {|\|\zeta_{y}\|-\|\zeta_x\||}{\|\zeta_{x}\|} \leq \frac{6\eps}{1-\delta}
\end{align*}
for $d(x,y)\leq R$. As $\delta$ is fixed, independent of $\eps$ we can achieve arbitrarily small variation. 

\ref{l2} $\iff$ \ref{Hilbert} $\iff$ \ref{Ozawakernel}\footnote{The equivalence of \ref{l2} and \ref{Ozawakernel} is also proved in \cite{Tu} Proposition 3.2.} $\iff$ \ref{Roekernel}: Note \ref{l2} $\implies$ \ref{Hilbert} trivially. To show \ref{Hilbert} $\implies$ \ref{Ozawakernel}, given $\xi_x$ let $u(x,y)=\Re\langle \xi_x,\xi_y\rangle$. Clearly this is of positive type and satisfies the appropriate vanishing conditions. The identity
$$\|\xi_x-\xi_{y}\|^2=2-2u(x,y)$$
shows that $(R,\sqrt{2\eps})$-variation for $\xi_x$ is equivalent to $(R,\eps)$-variation for $u$. To show that \ref{Ozawakernel} $\implies$ \ref{Roekernel}, note that bounded geometry, along with the support condition of \ref{Ozawakernel} implies that convolution with $u$ defines a bounded operator $\Op(u)$ on $l^2(X)$, specifically an element of the uniform Roe algebra $C^*_u(X)$.

Now we'll show that \ref{Roekernel} $\implies$ \ref{l2}. Let $u$ be a kernel with $(R,\eps)$-variation, and $\Op(u)$ the corresponding operator in $C^*_u(X)$. As $u$ is of positive type, $\Op(u)$ is a positive operator so it has a positive square root. We denote the corresponding kernel by $v$ and note that as $\Op(v)$ lies in $C^*_u(X)$ there is a self-adjoint kernel $w$ satisfying the support condition and the inequality
$$\|\Op(v)-\Op(w)\|<\min\left(\eps,\frac{\eps}{2(\|\Op(v)\|+\eps)}\right).$$
{}From this it follows that $\|\Op(w)\|\leq\|\Op(v)\|+\eps$ so we have
$$\|\Op(v)^2-\Op(w)^2\| \leq \|\Op(v)\| \|\Op(v)-\Op(w)\|+\|\Op(v)-\Op(w)\|\|\Op(w)\|<\eps.$$

Let $\zeta_x\in l^2(X)$ be the vector with entries $\zeta_x(z)=w(z,x)$. Now observe that
$$\langle\zeta_x,\zeta_{y}\rangle = \sum_z \overline{w(z,x)}w(z,y) = \sum_z w(x,z)w(z,y),$$
i.e.\ the kernel $\langle\zeta_x,\zeta_{y}\rangle$ consists of the matrix entries of the operator $\Op(w)^2$. Since $\Op(w)^2$ differs from $\Op(v)^2=\Op(u)$ by at most $\eps$ it follows that the kernel $\langle\zeta_x,\zeta_{y}\rangle$ differs from $u(x,y)$ entrywise by at most $\eps$. (In fact this condition is much weaker.) Hence the kernel $\langle\zeta_x,\zeta_{y}\rangle$ has $(R,3\eps)$-variation, so $\zeta_x$ has $(R,\sqrt{6\eps})$ variation. Finally replace $\zeta_x$ by $\eta_x=\zeta_x/\|\zeta_x\|$ and note that $\|\zeta_x\|^2=\langle\zeta_x,\zeta_x\rangle\geq 1-2\eps$. As in the proof of \ref{l2deltaweak} $\implies$ \ref{l2} we conclude that $\eta_x$ has $(R,2\sqrt{6\eps/(1-2\eps)})$-variation.
\end{proof}

\section{Morita invariance of $C^*_u(X)$}

The uniform Roe algebra is not a coarse invariant, as the following example illustrates. The uniform Roe algebra of a finite space is $M_n(\C)$ where $n$ is the cardinality of the space, however all finite spaces are coarsely equivalent. While the algebras $M_n(\C)$ are not all isomorphic, they are all  Morita equivalent. Note that for unital $C^*$-algebras, such as the uniform Roe algebra, two algebras $A$ and 
$B$ are Morita equivalent if and only if they are \emph{stably isomorphic} \cite[Thm.~7.6]{L}, which means that $A\otimes \K\cong B\otimes \K$, where $\K$ denotes the algebra 
of compact operators.

We prove the following. 
\begin{thm}
\label{moritainvariance}
If $X$ and $Y$ are uniformly discrete bounded geometry spaces, and $X$ is coarsely equivalent to $Y$, then $C^*_u(X)$ is Morita equivalent to $C^*_u(Y)$. 
\end{thm}

Before proving the theorem, we note the following corollary.

\begin{cor}\label{Cor5}
If $X$ and $Y$ are uniformly discrete bounded geometry spaces, and $X$ is coarsely equivalent to $Y$ then $C^*_u(X)$ is nuclear (resp.\ exact) if and only if $C^*_u(Y)$ is nuclear (resp.\ exact).
\end{cor}

\begin{proof}[Proof of Corollary \ref{Cor5}]
By the Theorem, if $X$ and $Y$ are coarsely equivalent, then the algebras $C_u^*(X)$ and 
$C_u^*(Y)$ are Morita equivalent. 
By Proposition 6.2 of \cite{KW}, if $A$ and $B$ are Morita equivalent $C^*$-algebras, then 
$A$ is nuclear if and only if $B$ is nuclear. 

In the case of exactness, we argue as follows. Since the algebras $C^*_u(X)$ and $C^*_u(Y)$ 
are unital and Morita equivalent, they are stably isomorphic. Since $\K$ is nuclear (and hence exact), if the algebra  $C_u^*(X)$ is exact, then the tensor product
$C^*_u(X)\otimes \K$ is also exact. Therefore $C^*_u(Y)\otimes \K$ is exact. Since a 
subalgebra of an exact $C^*$-algebra is also exact \cite[Prop.~2.6]{W}, the inclusion 
$C^*_u(Y) \hookrightarrow C^*_u(Y) \otimes \K$, defined with the help of a  projection 
in $\K$, shows that $C^*_u(Y)$ is exact, as required. 
 \end{proof}

We will now prove the theorem.

\begin{proof}[Proof of Theorem \ref{moritainvariance}]
First we will prove this under the assumption that the coarse equivalence $f\colon X \to Y$ is surjective. The issue is one of multiplicities, as in the above example of the finite spaces; if the map $f$ is actually a bijection, then we would have an isomorphism from $C^*_u(X)$ to $C^*_u(Y)$. We will prove that  $C^*_u(X)\otimes \K(l^2(\Z)) \cong C^*_u(Y)\otimes \K(l^2(\Z))$, where $\K(l^2(\Z))$ denotes the algebra of compact operators on the Hilbert space $l^2(\Z)$;  this is equivalent to Morita equivalence of the algebras $C^*_u(X)$ and $C^*_u(Y)$. Note that the algebras $C^*_u(X)\otimes \K(l^2(\Z))$ and $C^*_u(Y)\otimes \K(l^2(\Z))$ can be viewed as algebras of operators on $l^2(X\times \Z)$ and $l^2(Y\times \Z)$ respectively.

Since $f$ is a coarse equivalence there exists $R>0$ such that for each $y$ in $Y$, the preimage $f^{-1}(y)$ lies in some $R$-ball in $X$. By bounded geometry of $X$, there exists $N$ such that for each $y$, the cardinality of $f^{-1}(y)$ is at most $N$. Define $N(y)$ to be the cardinality of $f^{-1}(y)$, and for each $y$, enumerate the points of $f^{-1}(y)$, i.e. pick a bijection of $f^{-1}(y)$ with $\{1,\dots,N(y)\} \subseteq \{1,\dots,N\}$. We therefore obtain an identification of $X$ with a subset of $Y\times\{1,\dots,N\}$. Let $\pi$ denote the corresponding projection from $X$ to $\{1,\dots,N\}$, so that $x\mapsto (f(x),\pi(x))$ is the above identification.

We define a map $\phi$ from $X\times \Z$ to $Y\times \Z$, by $\phi(x,j)=(f(x),\pi(x)+jN(f(x)))$. Since for each $y$ in $Y$ there is exactly one $x$ in $X$ with $f(x)=y,\pi(x)=i$ for $i=1,\dots,N(y)$, the map $\phi$ is a bijection. This bijection gives rise to a unitary isomorphism from $l^2(X\times \Z)$ to $l^2(Y\times \Z)$, and hence an isomorphism $\Phi$ from the algebra of bounded operators on $l^2(X\times \Z)$ to the algebra of bounded operators on $l^2(Y\times \Z)$. We claim that $\Phi$ maps $C^*_u(X)\otimes \K(l^2(\Z))$ into $C^*_u(Y)\otimes \K(l^2(\Z))$, and $\Phi^{-1}$ maps $C^*_u(Y)\otimes \K(l^2(\Z))$ into $C^*_u(X)\otimes \K(l^2(\Z))$. Hence the restrictions of $\Phi$ to $C^*_u(X)\otimes \K(l^2(\Z))$ and $\Phi^{-1}$ to $C^*_u(Y)\otimes \K(l^2(\Z))$ give an isomorphism between $C^*_u(X)\otimes \K(l^2(\Z))$ and  $C^*_u(Y)\otimes \K(l^2(\Z))$ as required.

First we will show that $\Phi$ maps $C^*_u(X)\otimes \K(l^2(\Z))$ into $C^*_u(Y)\otimes \K(l^2(\Z))$. Since $\Phi$ is continuous, it suffices to prove this for the dense subalgebra of $C^*_u(X)\otimes \K(l^2(\Z))$ consisting of sums of elementary tensors of the form $T\otimes M$ where $T$ is a finite propagation operator on $l^2(X)$ and $M$ is a finite matrix. Indeed it is sufficient to prove that $\Phi(T\otimes e_{jj'})$ lies in $C^*_u(Y)\otimes \K(l^2(\Z))$, for $T$ of finite propagation, and where $e_{jj'}$ denotes a matrix unit. Partition $X$ as $X=\bigcup_{n=1,\dots,N,\,i=1,\dots,n} X_{n,i}$ where
$$X_{n,i}=\{x\in X : N(f(x))=n, \pi(x)=i\}.$$
We can write $T$ as the sum
$$T=\underset{i\leq n,i'\leq n'}{\sum_{n,n'\leq N}} P_{n,i}TP_{n',i'}$$
where $P_{n,i}$ denotes the projection of $l^2(X)$ onto the subspace $l^2(X_{n,i})$. Note that for each $i$, the restriction of $f$ to $X_i=\bigcup_{n\geq i}X_{n,i}$ is injective, and let $V_i$ denote the corresponding isometry from $l^2(X_i)$ to $l^2(Y)$. Fix $n,n'$ and $i,i'$, and let $S=P_{n,i}TP_{n',i'}$. Then $\Phi(S\otimes e_{jj'})=V_iSV_{i'}^*\otimes e_{i+nj,i'+n'j'}$. Since $f$ is a coarse equivalence, the operator $V_iSV_{i'}^*$ is of finite propagation, hence $V_iSV_{i'}^*\otimes e_{i+nj,i'+n'j'}$ lies in $C^*_u(Y)\otimes \K(l^2(\Z))$. Since this holds for each $n$, $n'$, $i$ and $i'$, we conclude that $\Phi(T\otimes e_{jj'})$ lies in $C^*_u(Y)\otimes \K(l^2(\Z))$ as required.

We will now show that $\Phi^{-1}$ maps $C^*_u(Y)\otimes \K(l^2(\Z))$ into $C^*_u(X)\otimes \K(l^2(\Z))$. As above, it suffices to show this for operators of the form $T\otimes e_{kk'}$, with $T$ a finite propagation operator on $l^2(Y)$. Let $Y_{n}=\{y\in Y : N(y)=n\}$, and let $P_n$ denote the projection of $l^2(Y)$ onto $l^2(Y_n)$. We can write $T$ as a sum $T=\sum_{n,n'\leq N}P_nTP_{n'}$. Now fix $n,n'$ and write $k=i+nj,k'=i'+n'j'$. Then for $S=P_nTP_{n'}$ we have $\Phi^{-1}(S\otimes e_{kk'}) = V_i^*SV_{i'} \otimes e_{jj'}$. As $f$ is a coarse equivalence $V_i^*SV_{i'}$ has finite propagation, and hence we conclude that $\Phi^{-1}$ maps $C^*_u(Y)\otimes \K(l^2(\Z))$ into $C^*_u(X)\otimes \K(l^2(\Z))$.

We have therefore shown that for $f\colon X \to Y$ a surjective coarse equivalence, $C^*_u(X)$ is Morita equivalent to $C^*_u(Y)$. The general case follows from the observation that given any coarse equivalence $f\colon X \to Y$, there are surjective coarse equivalences from both $X$ and $Y$ to the image $f(X)$. Hence we have the following isomorphisms:
$$C^*_u(X)\otimes \K(l^2(\Z)) \cong C^*_u(f(X))\otimes \K(l^2(\Z)) \cong C^*_u(Y)\otimes \K(l^2(\Z)).$$
\end{proof}

In \cite{Oz} Ozawa showed that a discrete group $G$ is exact if and only if its uniform Roe algebra 
is nuclear. Using Theorem \ref{moritainvariance} we can extend this result as follows.  
\begin{cor}
Let $G$ be a countable group acting properly by isometries on a uniformly discrete proper metric space $X$. If $C^*_u(X)$ is an exact algebra, then $G$ is an exact group. If moreover the action is cocompact, then the following are equivalent.
\begin{enumerate}
\item $X$ has property A;
\item $C^*_u(X)$ is nuclear;
\item $C^*_u(X)$ is exact;
\item $G$ is exact.
\end{enumerate}
\end{cor}

\begin{proof}
Pick a basepoint $x_0$ in $X$. Then $G$ is coarsely equivalent to the orbit $Y=Gx_0$. Exactness passes to subalgebras, so if $C^*_u(X)$ is exact, then $C^*_u(Y)$ is also exact. The algebra $C^*_u(Y)$ is Morita equivalent to $C^*_u(G)$, so the latter is exact. This algebra contains a subalgebra isomorphic to  $C^*_r(G)$, so $C^*_r(G)$ is also exact, hence $G$ is an exact group.

If the action of $G$ is cocompact then $G$ is coarsely equivalent to $X$, hence property A for $X$, nuclearity of $C^*_u(X)$ and exactness of $C^*_u(X)$ are all equivalent to their counterparts for $G$. For the group $G$, property A, nuclearity of $C^*_u(G)$, exactness of $C^*_u(G)$ and exactness of $G$ are all equivalent by the results of \cite{GK,HR,Oz}.
\end{proof}

\begin{rmk}
An alternative to the direct proof of Theorem \ref{moritainvariance} would be to rely on the following result from groupoid theory \cite[3.6]{STY}: If $X$ and $Y$ are coarsely equivalent uniformly locally 
finite coarse spaces then the groupoids $G(X)$ and $G(Y)$ are Morita equivalent. Given 
that $C^*_r(G(X))\cong C^*_u(X)$ one then appeals to the fact that Morita equivalence of 
groupoids implies the Morita equivalence of the associated reduced groupoid $C^*$-algebras
 \cite{MRW}. 
\end{rmk}

\section{Partial translation structures}\label{Section3}

A group $G$ comes equipped with both left and right multiplications. We assume
that a discrete group $G$ is equipped with a left invariant metric so that left multiplication acts by isometries on $G$.  
As noted in the introduction, the right multiplication by an element $g\in G$ moves every element 
$h$ of the group  the same distance.  We will reserve the term \emph{translation} for the right action on the group and the term  \emph{cotranslation} will refer to the left action.

We will now define analogues of both of these actions for a discrete metric space.

\begin{defn}
A \emph{partial bijection} from $X$ to $X$ is a subset $s$ of $X\times X$ such that the coordinate projections of $s$ onto $X$ are injective.
\end{defn}
A partial bijection can be viewed as a partially defined injection from $X$ into $X$, and we will write $x=s(y)$ if $(x,y)\in s$.

\begin{defn}
A \emph{partial translation} of $X$ is a partial bijection $t$ such that $d(x,y)$ is bounded for $(x,y)\in t$. The \emph{identity} translation, denoted $1$, is the diagonal of $X\times X$. The \emph{inverse} of $t$ is $t^*=\{(y,x):(x,y) \in t\}$.
\end{defn}
A partial translation gives rise to a partial isometry in $C^*_u(X)$; the inverse gives the adjoint partial isometry.

For $G$ a discrete group and for $g\in G$, the set $t_g=\{(h,hg)\colon h \in G\}$ is a (globally defined) partial translation. Viewing $t_g$ as a map from $G$ to itself, it acts by right multiplication by $g^{-1}$.

\begin{defn}
Let $\T$ be a collection of disjoint partial translations of $X$. A partial bijection $\sigma$ of $X$ is a \emph{partial cotranslation for $\T$} if for all $t\in\T$ and $(x,y)\in t$ such that $\sigma$ is defined on both $x$ and $y$, we have $(\sigma x,\sigma y) \in t$.
\end{defn}

For $G$ a group, elements of $G$ acting by left multiplication are (globally defined) partial cotranslations for $\{t_g : g\in G\}$.

\begin{defn}
\label{pt}
A \emph{ partial   translation structure} on $X$ is a collection $\T$ of partial translations of $X$, such that for all $R>0$ there is a finite subset $\T_R$ of disjoint partial translations in $\T$, and a collection $\SR$ of partial cotranslations of $\T_R$ satisfying the following axioms.
\begin{enumerate}
\item the union of the partial translations $t$ in $\T_R$ contains the $R$-neighbourhood of the diagonal, 
that is the set of all $(x,y)\in X\times X$ such that $d(x,y) < R$.

\item there exists $k$ such that for each $x,x'$ in $X$, there are at most $k$ elements $\sigma$ in $\SR$ such that $\sigma x=x'$;

\item for each $t$ in $\T_R$ and for all $(x,y),(x',y')$ in $t$, there exists $\sigma$ in $\SR$ such that $\sigma x=x'$ and $\sigma y=y'$.
\end{enumerate}
\end{defn}

A family $\{(T_R, \SR)\mid R>0\}$ of partial translations and partial cotranslations satisfying conditions 1-3 will be called an \emph{atlas}. Note that there is always a partial translation structure associated with an atlas given by putting $\T=\bigcup_{R>0}\T_R$. Indeed we can always enlarge this family $\T$ to include any additional partial translations we desire and the enlarged family will still be a partial translation structure. In particular if $X$ has any partial translation structure or equivalently any atlas, then taking $\T$ to be the family of all partial translations we obtain the maximal partial translation structure.

An atlas $\{(\T_R, \SR)\mid R>0\}$ where each $\T_R$ lies in a partial translation stucture $\T$ will be called an \emph{atlas for $\T$}. Note that atlases can be combined in the the following way: given any family $\mathcal A_i=\{(\T_R^i, \SR^i)\mid R>0\}$ of atlases for $\T$ we may construct a new atlas for $\T$ by choosing one of the pairs $(\T_R^i, \SR^i)$ for each $R>0$. 

For $\mathcal A=\{(T_R, \SR)\mid R>0\}$ an atlas on $X$, and for $R>0$, let $k_{\mathcal A}(R)$ denote the smallest $k$ such that for all $x,x'$ in $X$, there are at most $k$ elements $\sigma$ in $\SR$ with $\sigma x=x'$. For $\T$ a partial translation structure on $X$ let $k_{\T}(R)$ denote the minimum $k_{\mathcal A}(R)$ taken over all atlases $\mathcal A$ for $\T$. 

It is not immediately clear for which spaces partial translation structures exist. In Theorem \ref{translation-on-X} we will show that every  uniformly discrete bounded geometry metric space does admit a  partial   translation structure which allows us to make the following definition.

\begin{defn}
For $(X,d)$ a uniformly discrete bounded geometry  metric space let $\T$ denote the partial translation structure consisting of all partial translations on $X$. We define the \emph{translation invariant} of $X$ to be the function
$$\kappa_X(R) = k_\T(R).$$
\end{defn}

Note that if $d,d'$ are two coarsely equivalent metrics on $X$ then any partial translation structure for $(X,d)$ is a partial translation structure for $(X,d')$ and vice versa. Thus $\underset{R}\sup\,\kappa_X(R)$ is invariant under coarse equivalence of metrics.

We now introduce two properties for an atlas (freeness and global control) which we will need later.

\begin{defn}
An atlas $\mathcal A$ on $X$ is said to be \emph{free} if  $k_{\mathcal A}(R)=1$ for all $R>0$. 
\end{defn}

Since atlases, may be combined it follows that a partial translation structure $\T$ admits a free atlas if and only if $k_\T(R)=1$ for all $R>0$, and that $X$ admits a free atlas if and only if $\kappa_X(R)=1$ for all $R>0$.

\begin{defn}
An atlas $\{(T_R, \SR)\mid R>0\}$  is said to be \emph{globally controlled} if the partial cotranslation orbit
$$\{(x',y') : \text{ there exists $\sigma$ in $\SR$ such that $\sigma x=x',\sigma y=y'$} \}$$
is a partial translation for all $R>0$ and $x,y\in X$.
\end{defn}

Note that for $R>0$ and $x,y$ with $d(x,y)\leq R$ it is automatic that the orbit under $\SR$ is a partial translation; the cotranslation orbit of $(x,y)$ under $\SR$ is the unique $t$ in $\T_R$ such that $(x,y)\in t$.

\mbox{}

The following proposition describes the motivating example of a partial translation structure. 

 \begin{prop}\label{group-translation-structure}
 Let $G$ be a countable  discrete group. Then $G$ admits a canonical atlas which is free and 
 globally controlled. In particular,  $\kappa_G(R)=1$ for all $R$. 
\end{prop}

\begin{proof}
Since $G$ is a countable discrete group it may be equipped with a bounded geometry left invariant metric $d$. 
Let $t_g$ denote the partial translation $t_g= \{(x,xg): x\in G\}$ and for $R>0$ let $\T_R=\{t_g\mid d(e,g) < R\}$. Since the metric has bounded geometry $\T_R$ is a finite set. 

The pair $(x,y)$ belongs to  the $R$-neighbourhood of the diagonal if and only if $d(x,y) < R$, so $d(e, x^{-1}y) < R$ by the left invariance of the metric. But then $(x,y) \in t_{x^{-1}y}$ and 
$t_{x^{-1}y}\in \T_R$. Hence $\T_R$ satisfies condition 1 of  Definition  \ref{pt}. 

For each element $h\in G$ we have the bijection $\sigma_h:G\rightarrow G$ defined by $\sigma_h(x)=hx$. For each partial translation $t_{g}$ and each element $(x, xg)\in t_{g}$  we have  
$$(\sigma_h(x),\sigma_h(xg))=(hx, hxg) = (\sigma_h(x),\sigma_h(x)g)\in t_g.
$$ Thus the bijections $\sigma_h$ are partial cotranslations for $\T_R$.

For all $R>0$ let $\SR=\{\sigma_h : h\in G\}$. Clearly for any elements $x,x'\in G$ there is exactly one element of $\SR$, namely the cotranslation $\sigma_{x'x^{-1}}$, such that $\sigma_{x'x^{-1}}(x)=x'$ so the subsets $\SR$ satisfy condition 2 of Definition \ref{pt}. 

Let $t_g\in \T_R$ and $(x, xg), (x', x'g)\in t_g$. Then setting $h=x'x^{-1}$ we have $\sigma_{h}(x)=x'$ and $\sigma_{h}(xg)=x'g$ so condition 3 of Definition \ref{pt} is also fulfilled. Set 
$\A = \{(\T_R, \SR)\}$. 

 As noted above, there is exactly one element of $\SR$ such that $\sigma_h(x)=x'$, namely  $\sigma_{x'x^{-1}}$, so we have $k_{\A}(R)=1$ and so $\kappa_G(R)$ is also 1 for all $R$. This means that the 
 atlas is free. 

Finally if $x,y\in G$ then the cotranslation orbit $\{(\sigma_h(x),\sigma_h(y))\mid  \sigma_h\in \SR\}$ is equal to the set $t_{x^{-1}y}=\{(h, hx^{-1}y)\mid h \in G\}$, so in particular it is a partial translation. Thus the atlas is globally controlled.
\end{proof}

 We call the partial translation structure $\T=\bigcup_R\T_R$ the canonical partial translation structure on the group $G$.

We shall now demonstrate the existence of partial translation structures for arbitrary uniformly discrete 
bounded geometry metric spaces. We shall use the following terminology. 

\begin{defn}
A metric space $X$ is $R$-separated if for all distinct $x,y$ in $X$, we have $d(x,y)\geq R$.
\end{defn}

For metric spaces the following standard lemma will be needed to demonstrate the existence of partial translation structures.

\begin{lemma}
\label{bddgeom}
Let $X$ be a uniformly discrete bounded geometry metric space. Then for all $R$ the space $X$ can be written as a finite disjoint union of $R$-separated subsets.
\end{lemma}

\begin{proof}
This is a colouring argument. Given $R$ let $n$ be an upper bound on the cardinality of the $R$-balls of $X$. As a uniformly discrete bounded geometry metric space must be countable we can enumerate the points of $X$ as $x_1,x_2,\dots$. We will define inductively a colouring $c\colon X \to \{1,\dots,n\}$, such that $d(x,y)<R$ implies $c(x)\neq c(y)$. Let $c(x_1)=1$. Now suppose we have defined $c$ on the set $\{x_1,\dots,x_j\}$ in such a way that for $x,y$ in $\{x_1,\dots,x_j\}$ with $d(x,y)\leq R$, we have $c(x)\neq c(y)$. The set $B_R(x_{j+1})$ contains at most $n-1$ points from the set $\{x_1,\dots,x_j\}$, so there exists $i\in \{1,\dots,n\}$ with $c(x)\neq i$ for all $x$ in $\{x_1,\dots,x_j\} \cap B_R(x_{j+1})$. Define $c(x_{j+1})=i$. This extension also has the property that points $x,y$ with $d(x,y)\leq R$ are coloured differently, hence by induction we can extend $c$ to a colouring of $X$ such that $d(x,y)\leq R$ implies $c(x)\neq c(y)$. The sets $X_i=c^{-1}(\{i\})$ are $R$-separated as required.
\end{proof}

\begin{thm}\label{translation-on-X}
Let $X$ be a uniformly discrete bounded geometry metric space. Then $X$ admits a  partial   translation structure.
\end{thm}

\begin{proof}
Fix $R$ and write $X$ as a disjoint union, $X=X_1\cup\dots \cup X_n$, of $S$-separated sets, for $S>2R$. Let $t_{ij}$ be the set of pairs $(x,y) \in X_i\times X_j$ such that $d(x,y)\leq R$. For each $x\in X_i$ there is at most one $y\in X_j$ with $d(x,y)\leq R$, since $X_j$ is $S$ separated, with  $S>2R$. Conversely for each $y\in X_j$ there is at most one $x\in X_i$ with $d(x,y) \leq R$. Hence $t_{ij}$ is a partial translation. Note that $t_{ji}=t_{ij}^*$; $t_{ii}$ is the diagonal of $X_i\times X_i$; and the union of the sets $t_{ij}$ is the $R$-neighbourhood of the diagonal. Define $\T_R=\{t_{ij}: i,j=1,\dots,n\}$.

We will now define cotranslations. We will write $\SR$ as a union of sets $\SR^{ij}$, for $i,j=1,\dots,n$ with $i\leq j$. For each partial translation $t=t_{ij}$, pick a transitive permutation $\sigma_{ij}$ of the set $t$. This gives a partial cotranslation for $t$: for $(x,y)$ in $t$ we will  define $\sigma_{ij}$ on 
$X_i\cup X_j$. We define $\sigma_{ij} x$ and $\sigma_{ij} y$ such that $(\sigma_{ij} x,\sigma_{ij} y)=\sigma_{ij}(x,y)$. This is well defined as for each $x$ in $X_i$ there is at most one $y$ in $X_j$ with $(x,y)$ in $t$ and conversely for each $y$ there is at most one $x$. We define $\SR^{ij}$ to be the set of powers of these permutations $\{\sigma_{ij}^m\colon m\in \Z\}$. 

Note that $\sigma_{ij}^m$ is defined on pairs $(x,y)$ in $t_{ij},t_{ji},t_{ii},t_{jj}$ but not on any other $t$ in $\T_R$. For $t=t_{ij},t_{ji},t_{ii}$ or $t_{jj}$ and $(x,y)$ in $t$, it is clear that $(\sigma_{ij}^m x,\sigma_{ij}^m y)$ also lies in $t$, so the elements of $\SR=\bigcup_{i\leq j}\SR^{ij}$ are cotranslations for $\T_R$. For each $i\leq j$ and $x,x'\in X$ there is at most one cotranslation $\sigma\in \SR^{ij}$ such that $\sigma x=x'$, thus there are at most $k=n(n+1)/2$ elements of $\SR$ taking $x$ to $x'$. By construction, for each $i\leq j$, and $(x,y),(x',y')$ in $t_{ij}$ there exists $m$ such that $\sigma_{ij}^m x=x',\sigma_{ij}^m y=y'$, and the same holds for $t_{ji}=t_{ij}^*$. Thus $\T_R,\SR$ satisfy the required axioms. Now define $\T=\bigcup_R \T_R$.
\end{proof}

We conclude this section by discussing the question:  under what circumstances do there exist free or globally controlled atlases? We will show that uniform embeddings in groups give rise to these, and we will discuss the issue of which spaces are uniformly embeddable in groups. Throughout, the groups we consider will be countable discrete groups, equipped with a proper, left-invariant metric; such a metric is unique up to coarse equivalence.

\begin{thm}
\label{ptforgroups}
Let $X$ be a space admitting an injective uniform embedding into some discrete group $G$. Then $X$ admits a free and globally controlled atlas, and in particular $\kappa_X(R)=1$ for all $R>0$.
\end{thm}

\begin{proof}
Let $d_X$ denote the metric on $X$, and let $d_G$ denote the metric on $G$. Let $\phi\colon X \to G$ denote the injective uniform embedding of $X$ into $G$. We construct a partially defined action of $G$ on $X$ as follows. For $y$ in $X$ and $g$ in $G$ we define $g\act y=x$ if $\phi(y)g^{-1}=\phi(x)$, while $g\act y$ is undefined if $\phi(y)g^{-1}$ is not in the image of $\phi$. Note that $g\act y$ is uniquely determined if it exists, since $\phi$ is injective.

Fix $g\in G$. If $g\act y=x$, then
$$d_G(\phi(x),\phi(y))=d_G(\phi(y)g^{-1},\phi(y))=d_G(e,g).$$
As $\phi$ is a uniform embedding, there exists $R$ such that for $x,y$ with $d_G(\phi(x),\phi(y))=d_G(e,g)$ we have $d_X(x,y)\leq R$. Thus if $g\act y=x$ then $d_X(x,y)\leq R$. Hence for each $g\in G$ the action $g\act$, viewed as a partial bijection of $X$, defines a partial translation of $X$.

As $\phi$ is a uniform embedding, for all $R$ there exists $S$ such that if $d_X(x,y)\leq R$ then $d_G(\phi(x),\phi(y))\leq S$ i.e.\ $d_G(e,\phi(x)^{-1}\phi(y))\leq S$. Let $\T_R$ denote the set of partial translations $g\act$ with $d_G(e,g)\leq S$. These are disjoint partial translations whose union contains all $(x,y)$ with $d_X(x,y)\leq R$, i.e.\ the union contains the $R$-neighbourhood of the diagonal.

We will now define the cotranslations on $X$. We construct another partially defined action of $G$ on $X$. For $x$ in $X$ and $h$ in $G$ we define $h\lact x=x'$ if $h\phi(x)=\phi(x')$, while $h\lact x$ is undefined if $h\phi(x)$ is not in the image of $\phi$.  Let $\SR$ denote the set of partial bijections $h\lact$ for $h\in G$. We will show that these are partial cotranslation for $\T_R$. We must check that for each $g$ in $\T_R$, and $(x,y)$ in the corresponding partial translation $g\act$, if $(h\lact x,h\lact y)$ is 
defined, then it  also lies in the partial translation $g\act$. The pair $(x,y)$ lies in $g\act$ if and only if $x=g\act y$, and we then have $(h\lact x,h\lact y)=(h\lact (g\act y),h\lact y)$, where this is defined. Note that $h\lact(g\act y)=g\act(h\lact y)$ when both of these are defined, since the two actions arise from left and right multiplication in the group, so they commute. Moreover if $h\lact (g\act y),h\lact y$ are both defined then so is $g\act(h\lact y)$. Hence
$$(h\lact x, h\lact y)=(h\lact(g\act y),h\lact y)=(g\act(h\lact y),h\lact y)$$
whenever $(h\lact x, h\lact y)$ is defined, and $(g\act(h\lact y),h\lact y)$ lies in the partial translation $g\act$ as required.

We have shown that $h\lact$ acts on the partial translations, so it is a cotranslation. Moreover, the set of these cotranslations acts transitively on each partial translation in $\T_R$. To see this, note that for each partial translation $g\act$ in $\T_R$, and for any two pairs $(g\act y,y),(g\act y',y')$ in $g\act$, there exists $h$, namely $h=\phi(y')\phi(y)^{-1}$, such that $h\lact$ takes $(g\act y,y)$ to $(g\act y', y')$.

For each $x,x'$ we note that there is a unique cotranslation, namely $h\lact$ for $h=\phi(x')\phi(x)^{-1}$, such that $h\lact x$ equals $x'$. Thus $\{(\T_R, \SR)\mid R>0\}$ is a free atlas on $X$.

Finally, for any $x,y$ in $X$, if $x'=h\lact x$ and $y'=h\lact y$ then $\phi(x')^{-1}\phi(y')=\phi(x)^{-1}\phi(y)$. Let $R=d_G(e,\phi(x)^{-1}\phi(y))=d_G(e,\phi(x')^{-1}\phi(y'))$. As $\phi$ is a uniform embedding there exists $S$ such that $d(e,\phi(x')^{-1}\phi(y'))\leq R$ implies $d_X(x',y')\leq S$. Thus $d_X(x',y')$ is bounded for $(x',y')$ in the cotranslation orbit of $(x,y)$, so the atlas is globally controlled.
\end{proof}

\begin{cor}\label{cor20}
If $X$ is a uniformly discrete bounded geometry metric space for which $\kappa_X(R)>1$ for some $R$, then $X$ does not admit an injective uniform embedding into any countable group.
\end{cor}

\section{Uniform embeddings in groups}
Motivated by Theorem \ref{ptforgroups} we turn to the question of when a locally finite metric space
may admit a uniform embedding into a discrete group.  We shall use the following result from 
\cite[3.2]{DGLY} as our guiding principle: 
A locally finite metric space admits a uniform embedding in a Hilbert space 
if and only if it admits a local uniform embedding. 

In this section we  formulate a similar local to global principle for 
embeddings in groups and for this we need to define local variants of global properties of metric spaces. 

Given a countable, uniformly discrete metric space $X$, let $\Fin(X)$ be the disjoint union of all finite subsets of $X$. To be explicit about the metric, we can enumerate the finite sets as $X_1,X_2,\dots$, and view $\Fin(X)$ as the union of the sets $X_i\times\{i^2\}$ in $X\times \N$. The metric $d$  that $\Fin X$ inherits from the product metric on $X\times \N$ agrees with the given metric on each set $X_i$. 
If we then measure the distance between the subsets $X_i$  in the usual way then the metric has
the property that $d(X_i,X_j)$ tends to infinity as $i,j$ tend to infinity with $i\neq j$. Any other metric on $\Fin(X)$ with these properties will be coarsely equivalent to $d$.

We say that a metric space $X$ locally has property A if $\Fin X$ has property A. We say that $X$ is locally uniformly embeddable in a space $Y$ if $\Fin X$ is uniformly embeddable in $\Fin Y$. 
We note that this notion of local uniform embeddability in Hilbert space agrees with the definition of 
\cite{DGLY}.

As a first example of our local to global principle we have the following. 
\begin{prop}
A bounded geometry uniformly discrete metric space $X$ has property A if and only if it locally has property A.
\end{prop}

\begin{proof}
We use the characterisation of property A in terms of kernels, as in condition \ref{Ozawakernel} of Theorem \ref{propertyAcharacterisations}. Given a positive kernel $u$ on $X$ with $(R,\eps)$-variation, and vanishing for $d(x,y)>S$, we can produce a positive kernel on $\Fin X$ by defining kernels $u_i$ on $X_i$ to be the restrictions of $u$ to $X_i\times X_i$. Each individual $u_i$ has $(R,\eps)$-variation, however this need not be true for the collection taken together as a kernel on $\Fin X$. Let $I$ be the finite set of $i$ such that $X_i$ is within distance $R$ of some other $X_j$. We define a kernel $v$ on $\Fin X$ by 
$$
v(x,y)= \begin{cases} 1, & \text{if}\; x\in X_i,y\in X_j, i,j\in I ; \\
u_i(x,y), &  \text{if}\; x,y\in X_i, i\notin I; \\
0, & \text{otherwise.}
\end{cases} 
$$
It is easy to see that $v$ has $(R,\eps)$-variation, and vanishes for $d(x,y)$ sufficiently large. The kernel $v$ can be regarded as a block-diagonal matrix in which each block is positive. 
Thus property A implies local property A.

Conversely, given a positive kernel $v$ on $\Fin X$ with $(R,\eps)$-variation, and vanishing for $d(x,y)>S$, define $v_i$ to be the restriction of $v$ to the set $X_i$. View $v_i$ as a kernel on $X$ vanishing outside $X_i$. Choose a subsequence $i_j$ such that $X_{i_j}$ is an increasing sequence of sets whose union is $X$. By a diagonal argument we can extract a further subsequence such that the kernels converge pointwise. Let $u$ denote the limit kernel. As each $v_i$ is positive, and vanishes for $d(x,y)>S$ the same is true for $u$. To see that $u$ has $(R,\eps)$-variation, note that whenever $d(x,y)\leq R$ and $x,y\in X_i$ we have $|v_i(x,y)-1|<\eps$. By construction, for any fixed $x,y$, if $j$ is sufficiently large then $x,y\in X_{i_j}$, hence if $d(x,y)\leq R$ then $|v_{i_j}(x,y)-1|<\eps$ for $j$ sufficiently large. Thus $u$ has $(R,\eps)$-variation. Hence local property A implies property A.
\end{proof}

The following Theorem and its Corollary \ref{ptforlocalue} establish an analogous principle for 
embedding spaces in groups. Together they will provide a local version of Theorem \ref{ptforgroups}. 
\begin{thm}
\label{localue}
Let $X$ be a countable discrete metric space, and $G$ a countable group. Then $X$ is uniformly embeddable in $G$ if and only if $X$ is locally uniformly embeddable in $G$.
\end{thm}

\begin{proof}
To show that uniform embeddability implies local uniform embeddability is straightforward. The uniform embedding gives uniform embeddings of each finite subset $X_i$ of $X$ into some finite subset $G_i$ of $G$, uniformly in $i$. By replacing finite subsets $G_i$ of $G$ with larger finite subsets if necessary, we can arrange that $G_i\neq G_j$ for $i\neq j$. The disjoint union $\bigsqcup G_i$ may be equipped
with a metric in a similar way to the construction of the metric on $\Fin X$. By reindexing the sets $G_i$ 
we may extend this to an enumeration of the finite subsets of $G$. In this way we obtain a uniform
embedding of $\bigsqcup G_i$ into $\Fin G$, and so a uniform embedding of $\Fin X$ into  $\Fin G$. (Recall that the coarse structure on 
$\Fin (G)$ is independent of the enumeration chosen.)

The converse is a limiting argument. We will take $X_i$ to be an increasing sequence of finite subsets of $X$ whose union is $X$, so that $\bigsqcup_i X_i$ is a subset of $\Fin X$. By assumption we have a uniform embedding $\phi=\sqcup_i\phi_i$ of this into $\Fin G$.

Enumerate the pairs $(x,y)\in X\times X$,  as $(x_n,y_n)$ for $n=1,2,\dots$. As $\phi$ is a coarse map, for each $n$ there is a finite subset $F_n$ of $G$ such that if $x_n,y_n\in X_i$ then we have $\phi_i(x_n)^{-1}\phi_i(y_n)$ in $F_n$. As $F_1$ is finite, we can extract a subsequence of $i=1,2,\dots$ for which $\phi_i(x_1)^{-1}\phi_i(y_1)$ is constant. Similarly there are further subsequences for which this is constant for $(x_2,y_2),(x_3,y_3),\dots$. Thus by a diagonal argument there is a subsequence $\phi_{i_j}$ such that for each $n$ the sequence $\phi_{i_j}(x_n)^{-1}\phi_{i_j}(y_n), j=1,2,\dots$, is ultimately constant.

For $x,y\in X$ define $g_{xy}$ to be the limit of $\phi_{i_j}(x)^{-1}\phi_{i_j}(y)$. Note that as $$(\phi_{i_j}(x)^{-1}\phi_{i_j}(y))(\phi_{i_j}(y)^{-1}\phi_{i_j}(z))=\phi_{i_j}(x)^{-1}\phi_{i_j}(z)$$
for all $j$ we have the identity $g_{xy}g_{yz}=g_{xz}$. Similarly $g_{xx}=e$ for all $x$ and $g_{yx}=g_{xy}^{-1}$. As $\phi$ is a uniform embedding (i.e.\ each $\phi_i$ is a uniform embedding uniformly in $i$), for all $R$ there exists $S$ such that $d(x,y)\leq R$ implies $|\phi_i(x)^{-1}\phi_i(y)|\leq S$ for all $i$. Conversely $|\phi_i(x)^{-1}\phi_i(y)|\leq R$ implies $d(x,y)\leq S$ for all $i$, where $|\cdot|$ denotes the length function on the group. Thus for all $R$ there exists $S$ such that
$$d(x,y)\leq R \implies |g_{xy}|\leq S,$$
$$|g_{xy}|\leq R \implies d(x,y)\leq S.$$
Now fix a basepoint $x_0$ in $x$ and define $\psi\colon X\to G$ by $\psi(x)=g_{x_0x}$. Then $d(\psi(x),\psi(y))=|g_{x_0x}^{-1}g_{x_0y}|=|g_{xy}|$. The above inequalities thus show that $\psi$ is a uniform embedding.
\end{proof}

\begin{cor}
\label{ptforlocalue}
If a countable discrete space $X$ is locally uniformly embeddable in a countable group $G$ then $X$ contains a coarsely equivalent subset $Y$ which admits a free, globally controlled atlas.\end{cor}

\begin{proof}
If $X$ is locally uniformly embeddable in $G$ then $X$ is uniformly embeddable in $G$. Let $\phi$ denote the embedding. As this is a uniform embedding, there exists $R$ such that $d(x,y)\geq R$ implies $\phi(x)\neq \phi(y)$. Let $Y$ be a maximal $R$-separated subset of $X$. Then $Y$ is coarsely equivalent to $X$, and the restriction of $\phi$ to $Y$ is an injective uniform embedding, so by Theorem \ref{ptforgroups}, $Y$ admits a free, globally controlled atlas.
\end{proof}

We finish this section by constructing a universal space $\Gamma_{\u}$ with the property that 
$\Gamma_{\u}$ is uniformly embeddable in a countable discrete group $G$ if and only if each
uniformly discrete bounded geometry metric space $X$ uniformly embeds in some countable group $G_X$. 

We begin by showing that every discrete bounded geometry space can be uniformly embedded in a graph with bounded valences. Note that any bounded geometry space can be embedded in its total coarsening space, constructed in \cite{W}. This is a locally finite simplicial complex, hence any bounded geometry space embeds into a locally finite graph, viz.\ the 1-skeleton of this. The issue is to produce a graph with bounded valences. The following can be found in \cite[Prop.~5.1]{DGLY}. For the convenience of the reader we give here an alternative proof which emphasises the local nature of the question. 

\begin{prop}\label{embedding-in-graphs}
Let $X$ be a uniformly discrete bounded geometry metric space. Then there exists a graph $\Gamma$ for which each vertex lies in at most three edges, such that $X$ uniformly embeds in $\Gamma$.
\end{prop}
\begin{proof}
The idea is a telescope construction. Let $\Gamma_i$ denote the 1-skeleton of the $i$th Rips complex, i.e.\ $\Gamma_i$ is the graph whose vertices are the points of $X$ and for which there is an edge between $x$ and $y$ if and only if $d(x,y)\leq i$. We can then form a telescope graph by taking $\Gamma_0\sqcup\Gamma_1\sqcup\dots$ and for each $x\in X$ and $i=0,1,\dots$ adding an edge from the vertex $x$ in $\Gamma_i$ to the corresponding point in $\Gamma_{i+1}$. The inclusion of $X$ into this telescope as $\Gamma_0$ gives a uniform embedding, however in the telescope graph the number of edges emanating from a vertex is unbounded. To fix this we adjust the graphs $\Gamma_i$. Each vertex will be replaced by a linear graph. 

\begin{figure}
\begin{center}
  \includegraphics{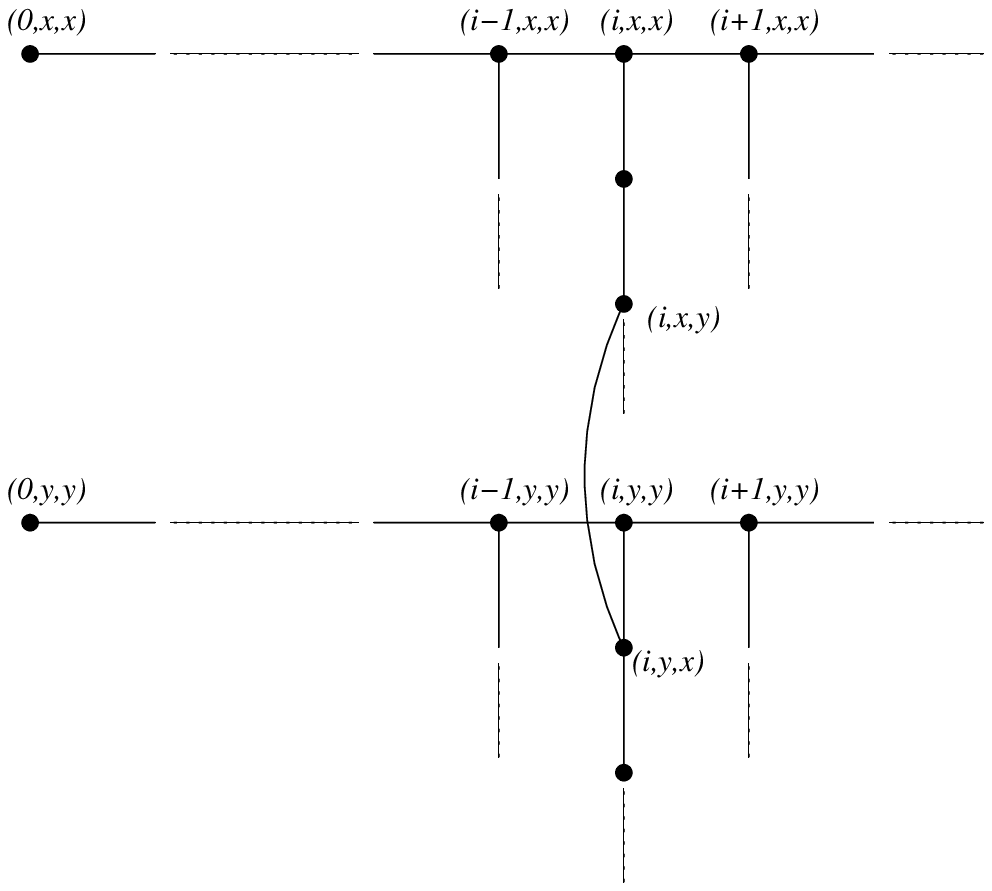}
\end{center}
\caption{The graph $\Gamma$.}
\end{figure}

We 
define the graph $\Gamma$ as follows. The vertex set of $\Gamma$ will be the set of triples $(i,x,y)$ with $i=0,1,2,\dots$, and $x,y$ in $X$ with $d(x,y)\leq i$. Let $b_{i,x}$ denote the set of triples $(i,x,y)$ with $d(x,y)\leq i$, and enumerate these starting with $(i,x,x)$. We use this enumeration to connect one
triple  to the next, thus $b_{i,x}$ is given the structure of a linear graph. Now for each $i$ and $x,y$ with $d(x,y)\leq i$, we add an edge joining $(i,x,y)$ with $(i,y,x)$ as Figure 1. If we fix $i$, and consider the projection $(i,x,y)\mapsto x$, this has the effect of collapsing each $b_{i,x}$ to a point, to recover the graph $\Gamma_i$. To complete the graph $\Gamma$, we carry out the telescope construction, by connecting $(i,x,x)$ to $(i+1,x,x)$ for each $i,x$. Note that the graph produced in this way has at most three edges emanating from any given vertex. Let $d_\Gamma$ denote the path length metric on $\Gamma$.

We will show that the inclusion of $X$ into $\Gamma$ given by mapping each $x\in X$ to the vertex  $(0,x,x)$ in $b_{0,x}$ is a uniform embedding. Denote this map by $\phi$. Given $R$, pick $i>R$, and note that as $X$ has bounded geometry, there is an upper bound $N$ on the number of points in an $i$-ball of $X$, thus, $|b_{i,x}|\leq N$ for all $x$ in $X$. Suppose that $x,y\in X$ with $d(x,y) \leq R$. Then there is a vertex $(i,x,y)$ in $b_{i,x}$ and a vertex $(i,y,x)$ in $b_{i,y}$, and these are connected by an edge. We have $d_\Gamma((i,x,x),(i,x,y))\leq N-1$, since $|b_{i,x}|\leq N$, and similarly $d_\Gamma((i,y,y),(i,y,x))\leq N-1$, thus $d_\Gamma((i,x,x),(i,y,y))\leq 2N-1$. On the other hand $\phi(x)=(0,x,x)$ and $\phi(y)=(0,y,y)$, so $d_\Gamma(\phi(x),(i,x,x))=i$ and $d_\Gamma(\phi(y),(i,y,y))=i$. Thus $d(x,y)\leq R$ implies $d_\Gamma(\phi(x),\phi(y))\leq S=2i+2N-1$.

Conversely, suppose that $d_\Gamma(\phi(x),\phi(y))\leq R$. Then there is a sequence of adjacent vertices $\phi(x)=v_0,v_1,\dots,v_k=\phi(y)$ in $\Gamma$ with $k\leq R$. Write $v_j$ in coordinates as $(i_j,x_j,y_j)$, and consider the sequence of points $x=x_0,x_1,\dots , x_k=y$ in $X$. Since $i_j\leq R/2$ for all $i$, we have $d(x_j,x_{j+1})\leq i_j \leq R/2$ for all $j$. It follows that if $d_\Gamma(\phi(x),\phi(y))\leq R$ then $d(x,y)\leq kR/2\leq S=R^2/2$. Thus $\phi$ is a uniform embedding.
\end{proof}

Let $\Gamma_\u$ denote the disjoint union of all connected finite graphs such that every vertex lies in at most three edges. This union may be equipped with a metric satisfying the following properties: 
1) the restriction to each graph is the standard edge metric; 2) for all $R> 0$ all but finitely many components  $\Gamma$ of  $\Gamma_\u$ have the property that the distance from $\Gamma$ 
to its complement in $\Gamma_\u$ is greater than $R$. This metric is unique  up to coarse equivalence. 

Note that $\Gamma_\u$ is a uniformly discrete metric space with   bounded geometry. Constructions of Gromov \cite{Gr} show that there exists a group $G$ and a coarse map from $\Gamma_\u$ to $G$, with a certain amount of control on how much the image of $\Gamma_\u$ is collapsed in $G$. This however is somewhat weaker than the assertion that $\Gamma_\u$ is uniformly embeddable in $G$. 
Proposition \ref{embedding-in-graphs} implies that $\Gamma_\u$ has the universal property that every bounded geometry space locally uniformly embeds in $\Gamma_\u$. 
Combining Proposition \ref{embedding-in-graphs} with Theorem \ref{localue} we 
obtain the following. 
\begin{thm}
\label{Gamma3}
The following are equivalent:
\begin{enumerate}
\item There exists a countable discrete group $G$ such that 
the space  $\Gamma_\u$ uniformly embeds in $G$; 
\item For each  uniformly discrete bounded geometry metric space $X$ there 
is a  countable discrete group $G_X$ such that $X$ uniformly embeds in $G_X$; 
 \item There exists a  countable discrete group $G$ such that every  uniformly discrete bounded geometry metric space $X$  uniformly embeds in $G$. 
 \end{enumerate}
 \end{thm}

\section{Translation algebras}\label{translation-algebras}
In this section we show how to associate a $C^*$-algebra $C^*(\T)$ to any partial translation structure $\T$ on a space $X$. This algebra will play the role assumed by the reduced $C^*$-algebra for a group.

\begin{defn}
For $\T$ a partial translation structure, the \emph{partial translation algebra} $C^*(\T)$ is the subalgebra of $C^*_u(X)$ generated by the partial translations $t\in \T$ (viewed as partial isometries).
\end{defn}

Note that given an atlas $\A=\{(\T_R,\SR)\mid R>0\}$ there is a canonical partial translation algebra $C^*(\T)$ associated to the partial translation structure $\T=\bigcup_R \T_R$. When $\mathcal A$ is free and globally controlled this algebra shares many of the important properties of the reduced $C^*$ algebra of a group. In particular it allows us to deduce properties of $C^*_u(X)$ from those of $C^*(\T)$, as in the proof of Theorem \ref{main} below.
 
Let $G$ be a countable discrete group and let $\T$ be the canonical partial translation structure on $G$ defined in Section \ref{Section3}. Note that $C^*(\T)$ is a subalgebra of $C^*_u(G)$, however $C^*_r(G)$ is not in general a subalgebra of $C^*_u(G)$. This is potentially a cause of confusion: The uniform Roe algebra for the group $G$ equipped with a \emph{right} invariant metric contains the algebra $C^*_r(G)$ which is the closure of the left regular representation. However, if $G$ is given the usual left invariant metric, then the uniform Roe algebra contains the right regular representation of $G$. This is because, with the left invariant metric, the left regular representation of the group acts by cotranslations, and not by operators of finite propagation, (i.e.\ translations) while the right regular representation of $G$ on $l^2(G)$ \emph{does} act by translations. It is this algebra which is recovered by the construction of  $C^*(\T)$.

\begin{thm}
Let $G$ be a countable discrete group and let $\T$ be the canonical partial translation structure on $G$. Then the algebra $C^*(\T)$ is canonically isomorphic to $C^*_r(G)$.
\end{thm}

\begin{proof}
It is immediate from the definitions that the algebra $C^*(\T)$ is equal to the closure of the right regular representation of $G$ on $l^2(G)$. The isomorphism of the closures of the left and right regular representations is given by conjugating by the unitary operator on $l^2(G)$ defined by $\delta_g\mapsto \delta_{g^{-1}}$ where $\delta_g$ denotes the usual basis vector corresponding to $g$.
\end{proof}

When the metric space
$X$ is sufficiently group-like, i.e., it admits a free and globally controlled atlas $\{(\T_R, \SR)\mid R>0\}$,  
we can use the  partial translation algebra corresponding to $\T=\bigcup_R\T_R$ to show that property A is equivalent to exactness  of $C^*_u(X)$. In particular, 
if the space is uniformly embeddable 
in a group we will show that $C^*_u(X)$ is nuclear if and only if it is exact. In the case where the atlas is free but not necessarily globally controlled we will show that the inclusion $C^*(\T) \hookrightarrow C^*_u(X)$ is a nuclear embedding if and only if  $X$ has property $A$. We begin with the following lemma.

\begin{lemma}
\label{control}
Let $H$ be a Hilbert space, and $A$ a $C^*$-subalgebra of $B(H)$. Suppose $\iota \colon A \hookrightarrow C^*_u(X)$ is a nuclear embedding. Then for any finite subset $E$ of $A$ and $\eps>0$ there exists a finite rank completely positive map $\theta\colon B(H) \to C^*_u(X)$ and $S>0$ with $\|\theta(a)-\iota(a)\|<\eps$ for $a\in E$, and with $\theta(a)$ of propagation at most $S$ for all $a\in B(H)$.
\end{lemma}

\begin{proof}
By nuclearity of the embedding there are completely positive maps $\phi\colon B(H) \to M_n$ and $\psi\colon M_n \to C^*_u(X)$ such that $\|\psi\circ\phi(a)-\iota(a)\|<\eps/2$ for all $a\in E$, see 
\cite[Def.~6.1.2]{RS}. There is a bijection between completely positive maps from $M_n$ to $C^*_u(X)$ and positive elements of $M_n(C^*_u(X))$, thus we can identify $\psi$ with a positive matrix $T$ in $M_n(C^*_u(X))$, namely $T_{ij} = \psi(e_{ij})$, where $e_{ij}$ is the standard matrix unit. 

We now approximate $T$ by a matrix whose entries are finite propagation operators.

Let $C$ be the maximum of the norms $\|\phi(a)\|_1$ for $a\in E$, where $\|\phi(a)\|_1$ denotes the $l^1$ norm, $\sum_{i,j}|\phi(a)_{ij}|$. Let $W$ be the square root of $T$. Let $W'$ be a finite propagation element of $M_n(C^*_u(X))$ such that $\|W'-W\|<\min(\eps,\frac{\eps}{4C(\|W\|+\eps)})$ and let $T'=(W')^2$. Then $T'$ is a positive matrix, and
$$\|T'-T\|=\|(W')^2-W^2\| \leq \|W'\| \|W'-W\|+\|W'-W\|\|W\|<\frac{\eps}{2C}.$$
Now let $\psi'\colon M_n \to C^*_u(X)$ be the completely positive map corresponding to $T'$, and define $\theta=\psi'\circ\phi$. Then for $a\in E$ we have
\begin{align*}
\|\theta(a)-\iota(a)\| &\leq \|\psi'\circ\phi(a)-\psi\circ\phi(a)\|+\|\psi\circ\phi(a)-\iota(a)\|\\
&=\|\sum_{i,j}\phi(a)_{ij}(T'-T)_{ij}\|+\|\psi\circ\phi(a)-\iota(a)\|\\
&\leq C\|T'-T\|+\eps/2<\eps.
\end{align*}

Let $S$ be the maximum of the propagations of the entries of $T'$. Then $\theta(a)$ has propagation at most $S$ for all $a\in B(H)$ as required.
\end{proof}

Now we are in a position to establish the following:

\begin{thm}\label{main}
Consider the following statements.
\begin{enumerate}
\item $X$ has property $A$.
\label{propA2}

\item $C^*_u(X)$ is nuclear.
\label{nuclear}

\item For every  partial   translation structure $\T$ on $X$, $C^*(\T)\hookrightarrow C^*_u(X)$ is a nuclear embedding.
\label{alltranslations}

\item There exists a partial translation structure $\T$ on $X$ admitting a free atlas and for which $C^*(\T)\hookrightarrow C^*_u(X)$ is a nuclear embedding.
\label{freetranslation}

\item $C^*_u(X)$ is exact.
\label{exact}

\item For every  partial   translation structure $\T$ on $X$, $C^*(\T)$ is exact.
\label{alltranslationsexact}

\item There exists a partial translation structure $\T$ on $X$ admitting an atlas which is free and globally controlled and for which $C^*(\T)$ is exact.
\label{globalfreetranslation}
\end{enumerate}
For $X$ a uniformly discrete  bounded geometry metric space we have (\ref{propA2})$\Longleftrightarrow$(\ref{nuclear})$\implies$(\ref{alltranslations}) and (\ref{nuclear})$\implies$(\ref{exact})$\implies$(\ref{alltranslationsexact}). If $\kappa_X(R)=1$ for all $R$, then conditions (\ref{propA2}), (\ref{nuclear}), (\ref{alltranslations}) and (\ref{freetranslation}) are all equivalent. If moreover there is a partial translation structure which admits a free and globally controlled atlas then (\ref{propA2}) through (\ref{globalfreetranslation}) are all equivalent.
\end{thm}

\begin{proof}
\ref{propA2} $\implies$ \ref{nuclear}: This is proved in \cite{R}, Proposition 11.41 and the equivalence $1\Leftrightarrow 2$ was proved in \cite{STY}. Here we give an alternative proof. We use the characterisation of property A, in terms of kernels (condition \ref{Ozawakernel} of Theorem \ref{propertyAcharacterisations}). For each $i$ there exists a positive type kernel $u_i$ such that 
\begin{enumerate}
\item $u_i(x,x)=1$; 
\item $|1-u_i(x,x')|\leq 1/i$, for $d(x,x') \leq i$;
\item  there exists $S_i$ such that $u_i(x,x')$ vanishes for $d(x,x')>S_i$.
\end{enumerate}

Define $\theta_i\colon \B(l^2(X))\to \B(l^2(X))$ to be the Schur multiplication by $u_i$. This is a completely positive contraction. The support condition on $u_i$ ensures that the range of $\theta_i$ lies in $C^*_u(X)$, while the variation condition ensures that for $T$ in $C^*_u(X)$ such that $d(x,x')\leq i$ on the support of $T$, we have $\|T-\theta_i(T)\|\leq \|T\|/i$. Thus for any $T$ in $C^*_u(X)$, $\theta_i(T)$ tends to $T$ as $i\to \infty$.

Now define $\Phi\colon \B(l^2(X)) \to l^\infty(X)$ to be restriction to the diagonal. Note that for a partial translation $t$, viewed as a partial isometry of $l^2(X)$, the operator $t\Phi(t^{*}T)$ is the restriction of $T$ to $t$. By bounded geometry, there exists a finite set $F_i$ of disjoint partial translations whose union includes all points with $d(x,x')\leq S_i$. For $t\in F_i$ define $\eta^t_i\colon \B(l^2(X)) \to l^\infty(X)$ by $\eta^t_i(T)=\Phi(t^{*}\theta_i(T))$, and note that this is a complete contraction. Then we have $\theta_i(T)=\sum_{t\in F_i}t\eta^t_i(T)$.

Now for any algebra $B$ we will show that the quotient map $Q\colon C^*_u(X)\otimes_{\mathrm{max}}B \to C^*_u(X)\otimes_{\mathrm{min}}B$ is injective. As $\eta^t_i$ is completely contractive, there is a well-defined map
$$\eta^t_i\otimes 1 \colon C^*_u(X)\otimes_{\mathrm{min}}B \to l^\infty(X)\otimes_{\mathrm{min}}B = l^\infty(X)\otimes_{\mathrm{max}}B \subset C^*_u(X)\otimes_{\mathrm{max}}B.$$
We use here the fact that $l^\infty(X)$ is nuclear. There is also a map
$$\theta_i\otimes 1 \colon C^*_u(X)\otimes_{\mathrm{max}}B \to  C^*_u(X)\otimes_{\mathrm{max}}B,$$
and we have $\theta_i\otimes 1=\sum_{t\in F_i}(t_i\eta^t_i\otimes 1)Q$; this holds for a dense subset of $C^*_u(X)\otimes_{\mathrm{max}}B$, and hence holds for all elements of $C^*_u(X)\otimes_{\mathrm{max}}B$.

If $Q(S)=0$ for $S \in C^*_u(X)\otimes_{\mathrm{max}}B$ then $(\theta_i\otimes 1)(S)=0$ for all $i$. But on the other hand we know that $(\theta_i\otimes 1)(S) \to S$ as $i\to \infty$, hence $Q(S)=0$ implies $S=0$. Injectivity of $Q$ implies nuclearity of $C^*_u(X)$.

\ref{nuclear} $\implies$ \ref{alltranslations} is immediate, as is \ref{nuclear} $\implies$ \ref{exact} $\implies$ \ref{alltranslationsexact}.

Now suppose $\kappa_X(R)=1$ for all $R$. Then there exists a partial translation structure with a free atlas, so \ref{alltranslations} $\implies$ \ref{freetranslation}.

\ref{freetranslation} $\implies$ \ref{propA2}: Again we use the characterisation of property A in terms of kernels. We must show that there is an approximate unit consisting of finite width positive kernels, for the algebra of functions on $X\times X$ tending to zero away from the diagonal. Nuclearity of the embedding provides approximate identity maps from $C^*(\T)$ into $C^*_u(X)$. Let $\delta_x, x\in X$ denote the standard basis for $l^2(X)$. The idea of the proof is to write the constant kernel 1 as $\langle \delta_x,t_{xy}\delta_y\rangle$ for certain operators $t_{xy}$, and then to use nuclearity to approximate the operators $t_{xy}$ in such a way as to produce kernels of finite width.

We make the following claim.
\begin{claim}
Let $X$ be a uniformly discrete bounded geometry metric space, and let $\T$ be a partial translation structure on $X$ with a free atlas $\{(\T_R,\SR)\mid R>0\}$. Then for all $R$, there exists a positive matrix $(t_{xy})$ indexed by $x,y$ in $X$ with entries in $C^*(\T)$, such that if $d(x,y)\leq R$ then $\langle \delta_x,t_{xy}\delta_y\rangle=1$, and such that the set of operators $\{t_{xy} : x,y\in X, d(x,y)\leq R\}$ is finite.
\end{claim}
Given this claim we can produce a positive kernel as follows. Fix $R,\eps>0$. As $C^*(\T)\hookrightarrow C^*_u(X)$ is a nuclear embedding we can apply Lemma \ref{control} to produce $S>0$, and a completely positive map $\theta$ from $B(l^2(X))$ into elements of $C^*_u(X)$ with propagation at most $S$, such that for $x,y$ with $d(x,y)\leq R$ we have $\|\theta(t_{xy})-t_{xy}\|<\eps$.

Define $u(x,y)=\langle \delta_x,\theta(t_{xy})\delta_y\rangle$. This is positive by positivity of the map $\theta$ and the matrix $(t_{xy})$. It vanishes if $d(x,y)>S$ as $\theta(t_{xy})$ has propagation at most $S$ for all $x,y$. If $d(x,y)\leq R$ then $\|\theta(t_{xy})-t_{xy}\|<\eps$ and $\langle \delta_x,t_{xy}\delta_y\rangle=1$, so $|1-u(x,y)|=|\langle \delta_x,(t_{xy}-\theta(t_{xy}))\delta_y\rangle|<\eps$. Hence we have produced a positive kernel with $(R,\eps)$-variation supported in the $S$-neighbourhood of the diagonal. As we can do this for each $R$ and $\eps$, the space $X$ has property A.

\bigskip

It remains to prove the Claim. Fix $R$. For each $x\in X$ we define an operator $s_x\colon l^2(X) \to l^2(\SR)$ as follows. For $x'\in X$ we define $s_x(\delta_{x'})$ to be $\delta_\sigma$ where $\sigma x=x'$, if such an element exists, and $s_x(\delta_{x'})=0$ otherwise. The element $\sigma$ is unique if it exists, since the atlas is free.

For each fixed $x,y$, $s_x^*s_y$ is an operator  on $l^2(X)$ thus it can be viewed as a matrix indexed by $X$. It has matrix entries
$$\langle \delta_{x'},s_x^*s_y \delta_{y'}\rangle=\langle s_x\delta_{x'},s_y \delta_{y'}\rangle, x',y'\in X$$
taking the value 1 if there exists $\sigma$ such that $\sigma x=x',\sigma y=y'$ and 0 otherwise. By hypothesis, if $d(x,y)\leq R$ then there exists such a partial cotranslation taking $(x,y)$ to $(x',y')$ if and only if $(x,y),(x',y')$ lie in the same element $t$ in $\T_R$. Thus if $d(x,y)\leq R$ then $s_x^*s_y$ is the unique partial translation $t$ in $\T_R$ such that $(x,y)\in t$.

We now define $t_{xy}=s_x^*s_y$. This is positive by construction. If $d(x,y)\leq R$ then $t_{xy}$ is a partial translation in the finite set $\T_R$, and moreover $(x,y)\in t_{xy}$ so $\langle \delta_x,t_{xy}\delta_y\rangle=1$. This establishes the claim.

\bigskip
Finally suppose that there exists a partial translation structure which admits a free and globally controlled atlas. Then each of \ref{nuclear}, \ref{alltranslations}, \ref{exact} and \ref{alltranslationsexact} immediately implies \ref{globalfreetranslation}. To see that \ref{propA2} and \ref{freetranslation} each imply \ref{globalfreetranslation} note that we have already established that \ref{propA2} is equivalent to \ref{nuclear} and to \ref{freetranslation} in this setting.

\ref{globalfreetranslation} $\implies$ \ref{propA2}: The proof is essentially the same as \ref{freetranslation} $\implies$ \ref{propA2}. Define $t_{xy}$ as before. We have the weaker hypothesis that $C^*(\T)$ is exact, i.e.\ $C^*(\T)\hookrightarrow B(l^2(X))$ is a nuclear embedding. We make use of \cite{Oz} Lemma 2. This implies that for all $R,\eps$ there exists a completely positive finite rank map $\theta \colon C^*(\T) \to B(l^2(X))$ such that:
\begin{enumerate}
\item $\theta$ has the form $\theta(\centerdot)=\sum_{i=1}^d \langle \delta_{a_i},\centerdot\delta_{b_i}\rangle T_i$, for some $a_i,b_i$ in $X$, and $T_i$ in $B(l^2(Y))$, and
\item for all $x,y$ with $d(x,y)\leq R$ we have $\|\theta(t_{xy})-t_{xy}\|<\eps.$
\end{enumerate}
Again we define $u(x,y)=\langle \delta_x,\theta(t_{xy})\delta_y\rangle$. Positivity and $(R,\eps)$-variation follow as before. Each $t_{xy}$ is a partial isometry whose support is a cotranslation orbit, and these supports partition $X\times X$. The global control condition implies that the supports are controlled, i.e.\ $t_{xy}$ has finite propagation for all $x,y$. Let $S$ be the maximum of the propagations of the operators $t_{a_ib_i}$. Then for $d(x,y)> S$ we have $\theta(t_{xy})=0$, since the support of $t_{xy}$ must be disjoint from the supports of the operators $t_{a_ib_i}$, and in particular it does not contain $(a_i,b_i)$. It follows that $u(x,y)$ vanishes when $d(x,y)>S$ as required.
\end{proof}

We conclude with the following corollaries.
\begin{cor}
If $X$ admits a (local) uniform embedding into a countable discrete group, then the following are equivalent: $C^*_u(X)$ is nuclear; $C^*_u(X)$ is exact; $C^*(\T)$ is exact for all partial translation structures $\T$ on $X$; $X$ has property $A$.
\end{cor}

\begin{proof}
By Corollary \ref{ptforlocalue}, if $X$ admits a local uniform embedding into a countable group then $X$ admits a free, globally controlled atlas. The result now follows immediately from the theorem.
\end{proof}

\begin{cor}
If the universal space $\Gamma_\u$ admits a (local) uniform embedding into a countable discrete group, then for every bounded geometry metric space $X$ the following are equivalent: $C^*_u(X)$ is nuclear; $C^*_u(X)$ is exact; $C^*(\T)$ is exact for all partial translation structures $\T$ on $X$; $X$ has property $A$.
\end{cor}

\begin{proof}
If the universal space $\Gamma_\u$ admits a local uniform embedding into a countable discrete group then by Theorem \ref{Gamma3}, so does every bounded geometry metric space. The result now follows from the previous corollary.
\end{proof}

\end{document}